\documentclass[a4paper,11pt]{amsart}
\usepackage{graphicx}

\newtheorem{theorem}{\quad Theorem}[section]

\newtheorem{proposition}{\quad Proposition}[section]

\newtheorem{conj}{\quad Conjecture}[section]

\renewcommand{\le}         {\leqslant}

\renewcommand{\ge}         {\geqslant}

\renewcommand{\tilde}         {\widetilde}
\renewcommand{\hat}         {\widehat}
\newcommand{\myint}         {\int\limits}
\newcommand{\myiint}         {\iint\limits}

\usepackage{amssymb}
\usepackage{dsfont,latexsym}

\usepackage{mathrsfs}
\renewcommand{\mathcal}{\mathscr}

\newcommand{\N}{\mathds{N}}

\newcommand{\R}{\mathds{R}}

\newcommand{\Per}{\,{\rm{Per}}}

\numberwithin{equation}{section}

\begin{document}
\title[A fractional framework for perimeters and 
phase transitions]{A fractional framework \\ for perimeters and 
phase transitions}

\author{Enrico Valdinoci}
\address{
Dipartimento di Matematica\\
Universit\`a degli Studi di Milano\\
Via Cesare Saldini 50\\
20133 Milano (Italy)
}
\email{enrico.valdinoci@unimi.it}
\thanks{I am greatly indebted to
Bego{\~n}a Barrios, Luis Caffarelli, 
Serena Dipierro, Alessio Figalli,
Giampiero Palatucci, Ovidiu Savin and Yannick Sire:
the results outlined in this note are the fruit of
the very pleasant and stimulating collaboration with them
and I profited enormously from 
the possibility of having
them as mentors and coworkers.
This work is supported by the {\it ERC} project~$\epsilon$ ({\it Elliptic
Pde's 
and Symmetry of Interfaces and Layers
for Odd Nonlinearities})
and the {\it FIRB} project A\&B ({\it Analysis
and Beyond}).}
\date{}

\begin{abstract}
We review some recent results on minimisers of a non-local perimeter
functional, in connection with some phase coexistence models
whose diffusion term is given by the fractional Laplacian.
\end{abstract}

\maketitle

\section{The fractional perimeter}

A notion of fractional perimeter was introduced in~\cite{CRS}.
To introduce it in a soft way, we consider a (measurable) 
set~$E\subseteq\R^n$ (say with~$n\ge2$) and a bounded, open\footnote{In
the sequel, for simplicity, the domain~$U$ will be
often implicitly assumed connected and with smooth boundary.} 
set~$U$
as in Figure~1.
\medskip

\begin{center} \label{inter}
\includegraphics[width=4in]{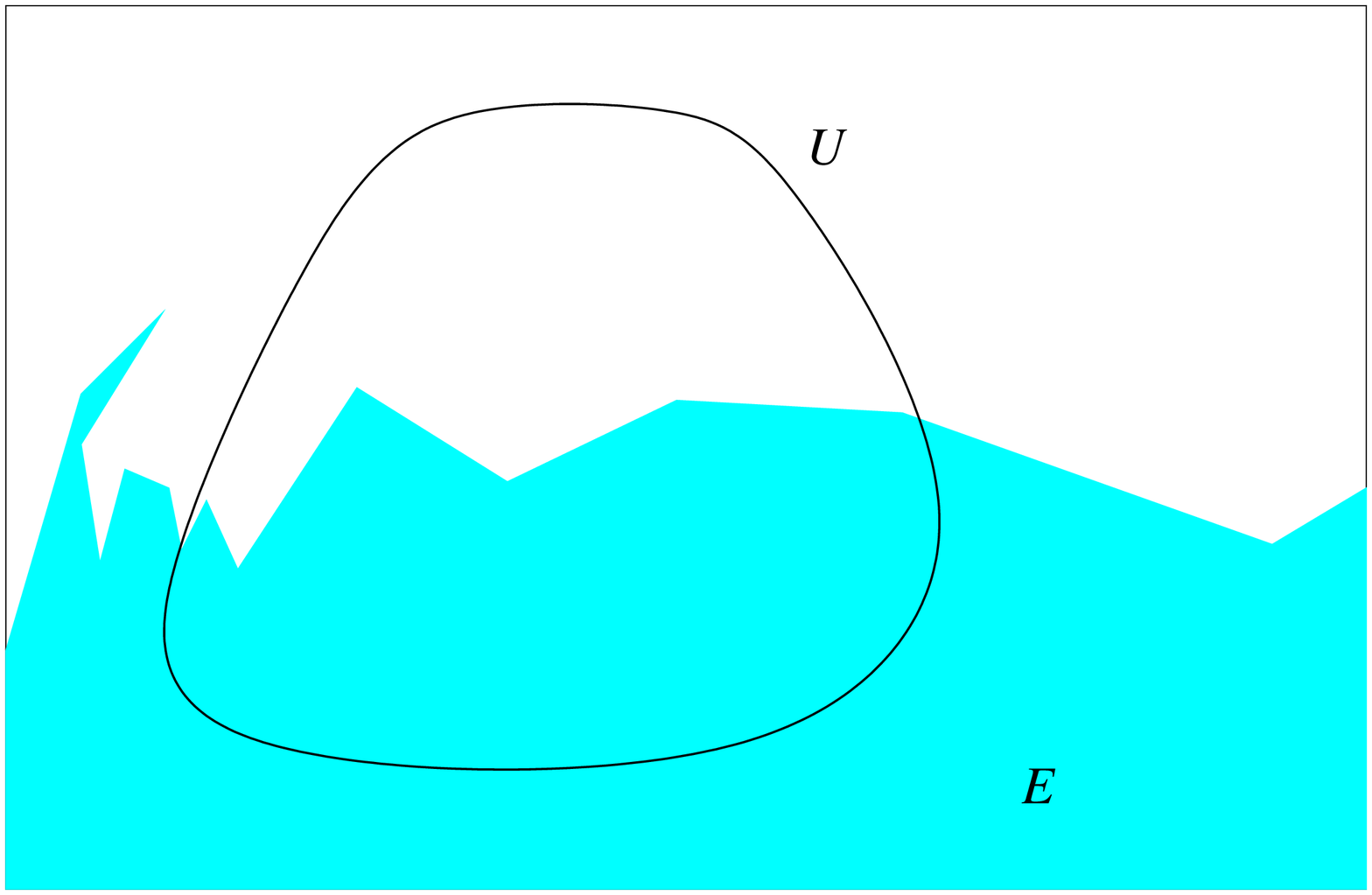}
\end{center}
\nopagebreak\centerline{\footnotesize\em
Figure~1. The sets~$E$ (in gray) and~$U$.}
\bigskip

The main idea of the fractional perimeter is that any point inside~$E$
``interacts'' with any point outside~$E$
giving rise to a functional whose minimisation
is taken into account. On the other hand, in the 
functional one may neglect the interactions that are fixed as ``boundary 
datum'' since they cannot contribute to the minimisation
(and they also may give an infinite contribution, which is
safer to take away). That is, the set~$O$ splits $E$ and its complement
into four sets, two inside~$E$, namely~$E':= E\cap U$ 
and~$E'':=E\setminus 
U$, and two outside~$E$, namely~$O':=
U\setminus E$ and~$O'':=(\R^n\setminus E)\cap (\R^n\setminus U)$,
see Figure~2. Then the functional
is the collection of the interactions of the points in~$E'$ and~$E''$
with the points in~$O'$ and~$O''$, with the exception of the 
interactions of points in~$E''$ with the ones in~$O''$, that are
``fixed by the boundary values''. 

\medskip

\begin{center} \label{inter2}
\includegraphics[width=4in]{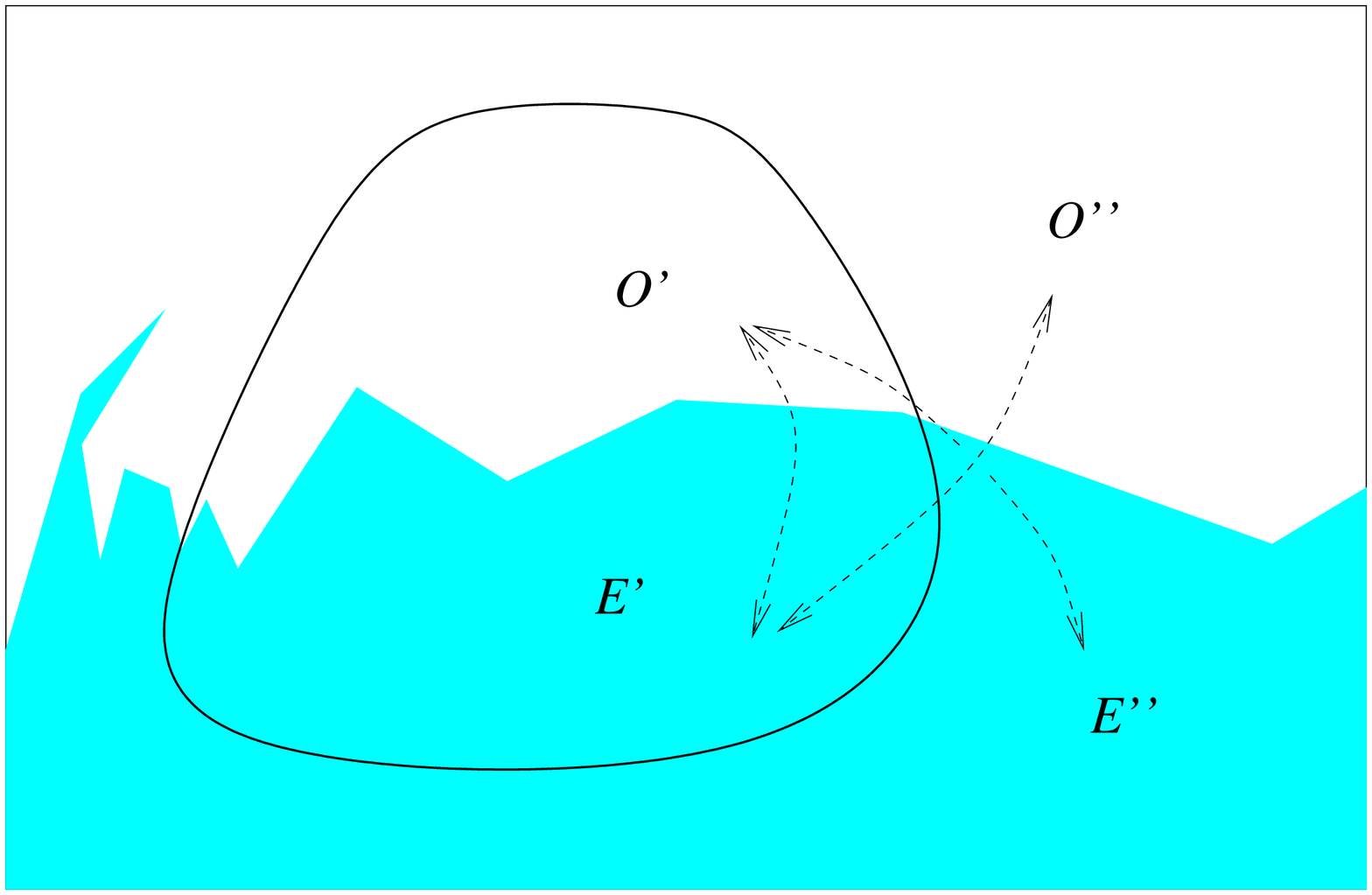}
\end{center}
\nopagebreak\centerline{\footnotesize\em
Figure~2.
The sets~$E'$, $E''$, $O'$ and~$O''$ and their interactions.}
\bigskip

Namely, one considers the functional
\begin{equation}\label{P}
\Per_s(E,U):= {\mathcal{I}}(E',O')+{\mathcal{I}}(E',O'')+
{\mathcal{I}}(E'',O'),\end{equation}
which formally coincides with~${\mathcal{I}}(E,\R^n\setminus E)-
{\mathcal{I}}(E'',O'')$ (though the latter may have no sense
since both~${\mathcal{I}}(E,\R^n\setminus E)$
and~${\mathcal{I}}(E'',O'')$ could be
infinite!). The interaction~${\mathcal{I}}$ that was considered
in~\cite{CRS} is
$$ {\mathcal{I}}(A,B):=\myiint_{A\times B}\frac{dx\,dy}{|x-y|^{n+2s}}$$
for any disjoint, measurable sets~$A$, $B$ and for a 
fixed~$s\in(0,1/2)$. 
The restriction on the range of~$s$ is natural, since the integrals
in~\eqref{P} diverge in general when~$s\in(-\infty,0]\cup[1/2,+\infty)$
(more precisely, for~$s\le0$ the contributions at infinity become
unbounded, while for~$s\ge1/2$ the problem arises from point~$x$ and~$y$
arbitrarily close to each other).
The functional in~\eqref{P} naturally produces
a minimisation problem: that is, one says that~$E$ is $s$-minimal in~$U$
if~$\Per_s(E,U)\le\Per_s(F,U)$ for any measurable set~$F$ that coincides
with~$E$ outside~$U$ (i.e., $F\setminus U=E\setminus U$).

The necessary compactness and semicontinuity properties
to ensure the existence of such $s$-minimisers are proved
in Section~3 of~\cite{CRS}, and the following result was obtained:

\begin{theorem}[Theorem 3.2 in \cite{CRS}]\label{99}
Let~$U\subset\R^n$
be a bounded Lipschitz domain and $E_o\subset\R^n\setminus U$
be a given set. There exists
a set~$E$, with~$E \setminus U
= E_o$ such that
$$ \Per_s(E,U)\le\Per_s(F,U)$$
for any~$F$ such that~$F\setminus U=E_o\setminus U$.
\end{theorem}

Moreover, in~\cite{CRS} $s$-minimisers are 
proved
to satisfy a suitable integral equation, that is the Euler-Lagrange
equation corresponding to the functional in~\eqref{P}. Namely
suppose that~$E$ is $s$-minimal in~$U$ and
that~$x_o\in U\cap (\partial E)$: then\footnote{We adopt
the standard notation for the characteristic function of a set~$E$, 
namely
$$\chi_E(x)=\left\{\begin{matrix}
1 & {\mbox{ if }}x\in E,\\
0 & {\mbox{ otherwise.}}
\end{matrix}\right.$$}
\begin{equation}\label{EL}
\myint_{\R^n}\frac{\chi_{E}(x_o+y)-
\chi_{\R^n\setminus E}(x_o+y)
}{|y|^{n+2s}}\,dy=0.\end{equation}
{F}rom the geometric point of view,~\eqref{EL} states that
a suitable average of~$E$ (centred
at any point of~$\partial E$)
is balanced by the average of its 
complement.
Due to the singularity of the denominator, \eqref{EL} only makes sense
for smooth sets: at this level, without knowing any a priori regularity
for the set~$E$, we must recall that equation~\eqref{EL}
must be taken in the viscosity sense (we refer to
Theorem~5.1 in~\cite{CRS} for details): in this setting, it may be 
interesting to notice that~\eqref{EL} says that~$(-\Delta)^s(\chi_{E}-
\chi_{\R^n\setminus E})=0$ along~$\partial E$ (see, e.g.,~\cite{guida}
for a basic introduction on the fractional Laplacian operator).

Of course, the functional in~\eqref{P}
may present a cumbersome combinatorics
which may complicate the computation of the interactions.
One may somehow turn around this difficulty by reducing the minimisation
problem in~\eqref{P} to a pde problem in~$\R^n\times(0,+\infty)$.
For this, given~$u:\R^n\rightarrow\R$, one introduces the extension 
of~$u$ as
$$ \tilde u(X):=\myint_{\R^n} {\mathcal{P}}(x-y,x_{n+1})u(y)\,dy\quad{\mbox{
with
}}\quad {\mathcal{P}}(X):=\frac{ c_{n,s} x_{n+1}^{2s}}{|X|^{n+2s}}.$$
Here we used the notation~$X:=(x,x_{n+1})\in\R^n\times(0,+\infty)$,
and~$c_{n,s}>0$ is a normalising constant. 
Given~$\Omega\subset\R^{n+1}$ and~$v:\R^n\times(0,+\infty)
\rightarrow\R$, we define
$$ {\mathcal{E}}_\Omega (v):=\myint_{\Omega\cap \{x_{n+1}>0\}}
x_{n+1}^{1-2s}|\nabla v(X)|^2 \,dX.$$
Then we have the following result:

\begin{proposition}[Proposition 7.3 in~\cite{CRS}]\label{73}
The set $E$ is an $s$-minimiser 
in a ball~$B$ if and only if the 
extension~$\tilde u$ of $\chi_E-\chi_{\R^n\setminus E}$ satisfies
$$ {\mathcal{E}}_\Omega(\tilde u)\le
{\mathcal{E}}_\Omega(v)$$
for all bounded Lipschitz domains $\Omega\subset\R^{n+1}$
with~$\Omega\cap \{x_{n+1}=0\} \Subset B$ 
and all functions~$v$ that equal~$\tilde u$ in a neighbourhood 
of~$\partial\Omega$
and take the values~$\pm1$ on~$\Omega\cap \{x_{n+1}=0\}$. 
\end{proposition}

Notice that the restriction for~$v$ to take values~$\pm1$
``on the trace''~$\{x_{n+1}=0\}$
(that is, to ``agree with a set'' on~$\R^n$)
causes several analytical difficulties in the choice
of the admissible perturbations of~$\tilde u$.

One of the main questions addressed in~\cite{CRS}
is the one of the regularity of the $s$-minimisers:
that is, it is shown there that~$s$-minimisers are smooth sets
outside a singular set of~$(n-2)$-Hausdorff dimension:

\begin{theorem}[Main Theorem 2.4 in~\cite{CRS}]\label{s989}
If $E$ is an $s$-minimiser in~$B_1$, then $\partial E \cap B_{1/2}$ is, 
to the 
possible 
exception of 
a closed set~$\Sigma$
of finite~$(n-2)$-Hausdorff dimension, a $C^{1,\alpha}$-hypersurface 
around each 
of its 
points.\end{theorem}

Notice that one expects~$\partial E$ to be ``an $(n-1)$-dimensional
object'', so the statement that~$\Sigma$ is  ``an $(n-2)$-dimensional
object'' states that~$\Sigma$ is somehow negligible inside~$\partial E$,
hence~$\partial E$ is smooth near ``the majority of its points''.
On the other hand, Theorem~\ref{s989}
leaves many questions open: for instance, is there any singular set at
all? are there any singularities if the dimension~$n$ is small enough?
what do the $s$-minimisers look like?

Some of these questions will be discussed in further detail
in \S~\ref{8855}-\ref{88811}. As for the latter problem,
it is quite embarrassing to admit that, at the moment, there
is a real lack of explicit examples: indeed, the
only explicit $s$-minimiser known is the half-plane
(which is in fact a minimiser in any domain~$U$):
this was proved in Corollary~5.3 
in~\cite{CRS} and the proof is based on a comparison principle
(i.e., if an~$s$-minimiser is contained in some strip outside~$U$,
then it is contained in the same strip inside~$U$ too).
Notice that this lack of explicit examples
does not prevent $s$-minimisers from existing (recall
Theorem~\ref{99}).
In any case, at the moment, no $s$-minimiser with a non-void
singular set is known.
\medskip

\begin{center} 
\includegraphics[width=4in]{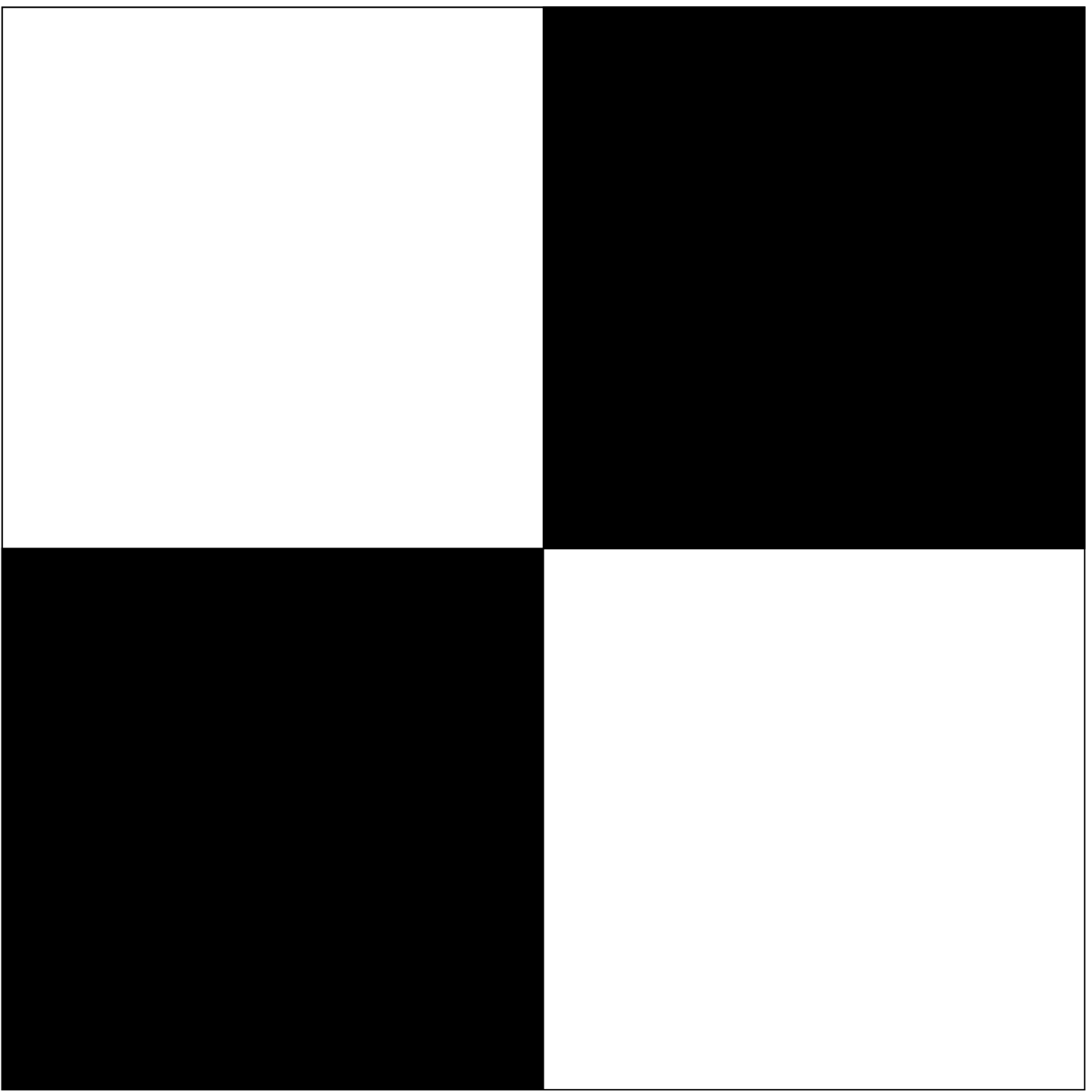}
\end{center}
\nopagebreak\centerline{\footnotesize\em
Figure~3.
The cone ${\mathcal{K}}$.}
\bigskip

One may also wonder if there are sets satisfying the Euler-Lagrange
equation in~\eqref{EL} that possess a non-void singular set:
the answer is in the affirmative, and
a simple example is given by the classical cone in the plane
\begin{equation}\label{cone.1}
{\mathcal{K}}:=\{ (x,y)\in \R^2 {\mbox{ s.t. }} 
xy>0\},\end{equation}
see Figure~3. Of course, $ {\mathcal{K}}$ has a singularity
at the origin, and, by symmetry, one sees that
\begin{equation}\label{K el}
{\mbox{${\mathcal{K}}$
satisfies \eqref{EL} (possibly in the viscosity sense)}}.\end{equation}

\subsection{Asymptotics of the $s$-perimeter}

Up to now, the reason for which we think that the functional 
in~\eqref{P} is a ``fractional perimeter'' may seem mysterious
to a reader not familiar with the subject.
The motivation arises for the asymptotics as~$s\nearrow1/2$
in which the functional~$\Per_s$
(suitably renormalised) approaches the classical
perimeter (as usual, we use the notation~$\omega_{n-1}:= {\mathcal{H}}^{n-1}(S^{n-1})$
for the surface of the $(n-1)$-dimensional sphere):

\begin{theorem}[\cite{CV1, AdPM}]\label{888990}
$\,$
\begin{itemize}
\item Let~$\alpha\in(0,1)$, $R>0$, $s_k\nearrow1/2$
and~$E$ be a set with~$C^{1,\alpha}$-boundary in~$B_R$.
Then
$$ \lim_{k\nearrow+\infty}(1-2s_k)\Per_{s_k}(E,B_r)=\omega_{n-1} \Per(E,B_r)
\qquad{\mbox{a.e. }}r\in(0,R).$$
\item Let~$R>r>0$, $s_k\nearrow1/2$ and~$E_k$ be such that
$$ \sup_{k\in\N}(1-2s_k)\Per_{s_k}(E_k,B_R)<+\infty.$$
Then, up to subsequence,~$\chi_{E_k}$ converges in~$L^1(B_r)$
to~$\chi_E$, for a suitable~$E$ with 
finite perimeter in~$B_r$.
\item Let~$R>r>0$.
Let~$s_k\nearrow1/2$ and~$E_k$ be $s_k$-minimisers 
in~$B_R$,
with~$\chi_{E_k}$ converging in~$L^1(B_R)$
to~$\chi_E$. Then $E$ has minimal perimeter in~$B_r$.
Also, $E_k$ approach~$E$ uniformly in~$B_r$, meaning
that for any~$\epsilon>0$ there exists~$k_o$ (possibly depending on~$r$
and~$\epsilon$) such that if~$k\ge k_o$ then~$E_k\cap B_r$ 
and~$B_r\setminus E_k$ are contained, respectively, in an 
$\epsilon$-neighbourhood of~$E$ and of~$\R^n\setminus E$.
\end{itemize}
\end{theorem}

Also, one can show that the convergence of the functional in~\eqref{P}
to the classical perimeter as~$s\nearrow1/2$ holds in a 
suitable~$\Gamma$-convergence sense: see~\cite{AdPM}. In any case,
we hope that this motivates the notation of fractional perimeter
introduced in~\eqref{P}. On the other hand, when~$s\searrow0$,
it is conceivable that the functional in~\eqref{P} must approach,
in some sense, the Lebesgue measure~${\mathcal{L}}^n$ (up to scaling).
To see this, let us recall the notion of Gagliardo seminorm
of a function~$u$:
$$ 
[u]_{G,s}:=\sqrt{\myiint_{\R^{2n}}\frac{|u(x)-u(y)|^2}{|x-y|^{n+2s}}\,dx\,dy}.$$
By taking the Fourier transform, one sees that
$$ [u]_{G,s}^2=c(n,s)\myint_{\R^n}|\xi|^{2s} |\hat u(\xi)|^2\,d\xi,$$
for any~$u$ in the
Schwartz space $C^\infty_\downarrow(\R^n)$
of rapidly decreasing smooth functions: here above~$\hat u$
is the Fourier transform of~$u$ and~$c(n,s)$ is a suitable
normalising constant with the property that
$$ \lim_{s\searrow0} c(n,s)\,s =c_n$$
for an appropriate~$c_n>0$ (see, e.g., Proposition~3.4
and Corollary~4.2 in~\cite{guida}).
Therefore
\begin{equation}\label{97}\begin{split}&
\lim_{s\searrow0}s [u]_{G,s}^2=\lim_{s\searrow0}
c(n,s)\,s\myint_{\R^n}|\xi|^{2s} |\hat 
u(\xi)|^2\,d\xi\\ &\qquad=c_n \myint_{\R^n}|\xi|^{0} |\hat
u(\xi)|^2\,d\xi=
c_n \|\hat u\|_{L^2(\R^n)}^2=
c_n \| u\|_{L^2(\R^n)}^2,\end{split}\end{equation}
thanks to Plancherel Theorem. Though this formula
is obtained here for~$u\in C^\infty_\downarrow(\R^n)$, it holds true also
for functions~$u\in L^2(\R^n)$ for which~$[u]_{G,s_o}$ is finite
for some~$s_o\in(0,1)$ (see, e.g.,~\cite{Maz}
for a general theory in~$L^p$-spaces). In particular,
we may take~$u:=\chi_E$ in~\eqref{97} for a smooth~$E\subset U$ 
(the smoothness of~$E$ ensures that~$[\chi_E]_{G,s_o}$ is finite
and the boundedness of~$E$ that~$\chi_E\in L^2(\R^n)$): we conclude that
\begin{equation}\label{9066}\begin{split}
&\lim_{s\searrow0}
2s\Per_s(E,U)=
\lim_{s\searrow0}
2s\myiint_{E\times(\R^n\setminus E)}\frac{dx\,dy}{|x-y|^{n+2s}}=
\lim_{s\searrow0}s [\chi_E]_{G,s}^2
\\&\qquad=c_n \| \chi_E\|_{L^2(\R^n)}^2=c_n {\mathcal{L}}^n(E).\end{split}
\end{equation}
The asymptotic behaviour as~$s\searrow0$ in the general case
is slightly more complicated and it is dealt with in~\cite{DFPV}:
the main difficulties are that the limit may not even exist
and, when it does exist, it is a suitable convex combination of
the normalised Lebesgue measure of~$E\cap U$ with the one of~$U\setminus E$,
with the convex interpolation parameter depending on the ``shape of~$E$
outside~$U$''. More precisely, one introduces the parameter
\begin{equation}\label{a}
a(E):=\lim_{s\searrow0}\frac{2s}{\omega_{n-1}}\myint_{E\setminus 
B_1}\frac{dy}{|y|^{n+2s}}\end{equation}
and the normalised Lebesgue measure~${\mathcal{M}}(E):=\omega_{n-1}
{\mathcal{L}}(E)$.
Notice that
$$ a(E)\le\lim_{s\searrow0}\frac{2s}{\omega_{n-1}}\myint_{
\R^n\setminus
B_1}\frac{dy}{|y|^{n+2s}}=
\lim_{s\searrow0}2s
\myint_{1}^{+\infty}\frac{\rho^{n-1}\,d\rho}{\rho^{n+2s}}=1$$
hence~$a(E)\in[0,1]$.
Then, in some sense (that we will make
precise below in Theorem~\ref{ts 0}),
the asymptotic behaviour as~$s\searrow0$ is 
given by 
the formula
\begin{equation}\label{s0}
\lim_{s\searrow0}2s \Per_s(E,U)=\big(1-a(E)\big) {\mathcal{M}}(E\cap U)+
a(E){\mathcal{M}}(U\setminus E).
\end{equation}
Notice that when~$E$ is a smooth subset of~$U$, then~$E\subset B_R$
for some~$R>0$, hence
$$ a(E)\le
\lim_{s\searrow0}\frac{2s}{\omega_{n-1}}\myint_{
B_R\setminus
B_1}\frac{dy}{|y|^{n+2s}}=
\lim_{s\searrow0} 2s\myint_{1}^{R}\frac{\rho^{n-1}\,d\rho}{\rho^{n+2s}}=0,$$
i.e.~$a(E)=0$
and so~\eqref{s0}
boils down to~\eqref{9066}. Moreover, \eqref{s0} states that the $s$-perimeter
(suitably normalised, which has a non-local nature)
approaches as~$s\searrow0$ the convex combinations of two measures ``localised''
in~$U$ (namely~${\mathcal{M}}(E\cap U)$ and~${\mathcal{M}}(U\setminus E)$),
but the combination parameter~$a(E)$ takes into account the contribution
of~$E$ ``coming from infinity''. Though this is rather attractive,
a rigorous statement has to take into account the possibilities
that the above limits do not exist, and the precise result
on the asymptotics as~$s\searrow0$ reads as follows:

\begin{theorem}[Theorems 2.5 and 2.7 in~\cite{DFPV}]\label{ts 0}
Let~$E$ be such that~$\Per_{s_o}(E,U)<+\infty$ for some~$s_o\in(0,1/2)$
and suppose that the limit defining~$a(E)$ in~\eqref{a} exists.
Then the limit in~\eqref{s0} holds true.

Also, if~$\Per_{s_o}(E,U)<+\infty$ for some~$s_o\in(0,1/2)$
and~${\mathcal{L}}^n(E\cap U)\ne{\mathcal{L}}^n(U\setminus E)$,
then the existence of the
limit defining~$a(E)$ in~\eqref{a} is equivalent to~\eqref{s0}.
\end{theorem}

The existence condition on the limit
defining~$a(E)$ in~\eqref{a} cannot be removed from Theorem~\ref{ts 0},
since \cite{DFPV} also provides an example when such limit does not exist
(and the limit in~\eqref{s0} does not exist as well).
\medskip

\begin{center} 
\includegraphics[width=4in]{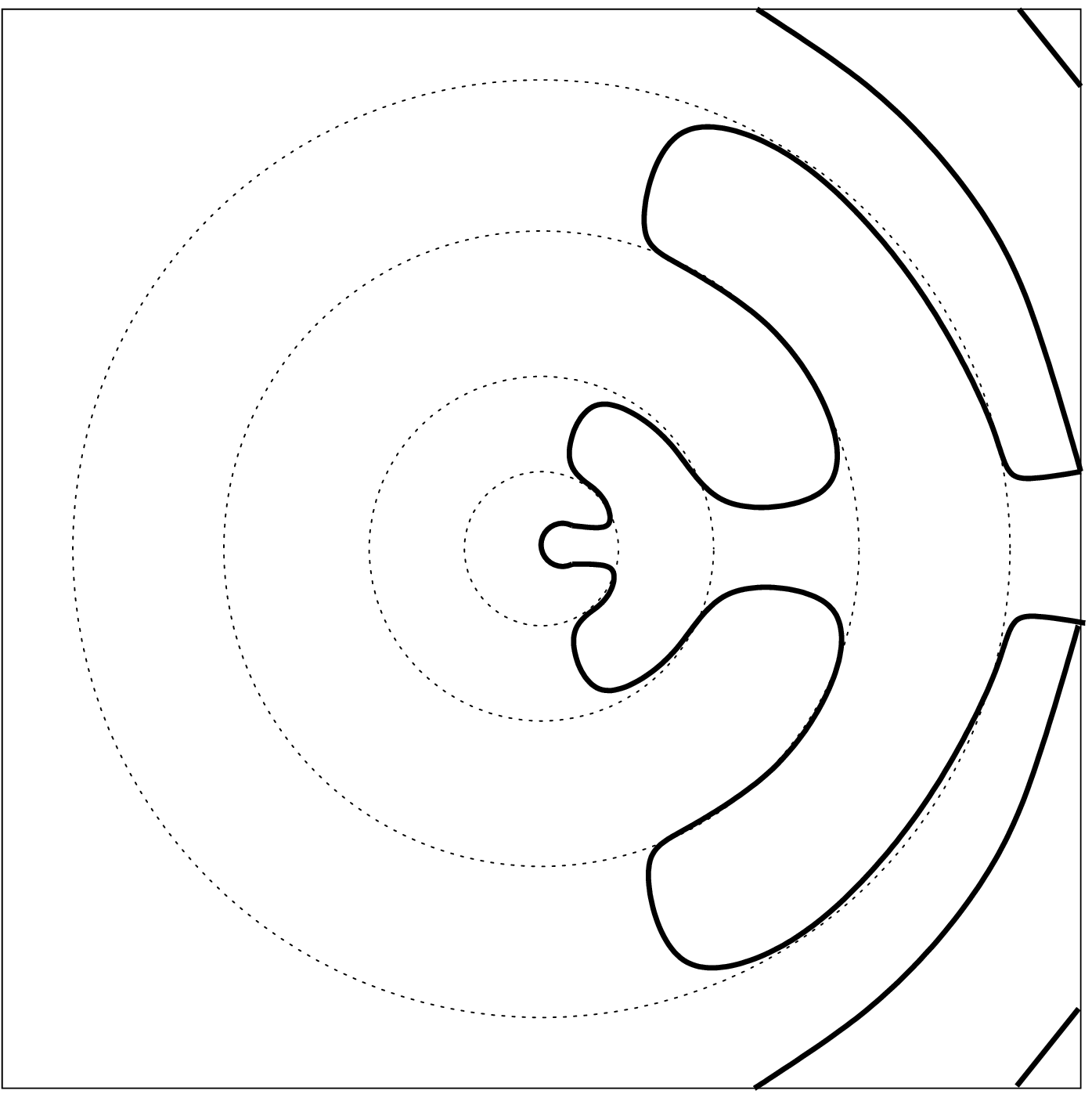}
\end{center}
\nopagebreak\centerline{\footnotesize\em
Figure~4.
An example for which the limit as~$s\searrow0$ of the fractional
perimeter does not exist.}
\bigskip

Very roughly speaking, the example (as grossly depicted in Figure~4)
considers 
a set~$E$
which ``looks like a cone'' of a small opening in the annulus~$B_{R_2}\setminus 
B_{R_1}$, then like a cone of a big opening in the annulus~$B_{R_3}\setminus 
B_{R_2}$, and so on, alternating cones of small and big openings in
subsequent annuli, with~$R_k\rightarrow+\infty$ to be chosen appropriately.
Then, the idea is that~$a(E)$ would ``feel alternatively'' the small
and the big cone openings in the asymptotics and consequently
the limit in~\eqref{a}
does not exist (of course, some computation is needed
to relate the ``spatial oscillation'' in the annuli with
the parameter~$s\searrow0$, see
Example~2.7 in~\cite{DFPV}
for details).
Anyway, from Theorems~\ref{888990}
and~\ref{ts 0}, with a slight abuse of notation, one may think that the
fractional perimeter interpolates the classical perimeter with
a weighted Lebesgue measure when the parameter $s$ varies in the 
range~$(0,1/2)$.

\subsection{Regularity of $s$-minimal sets in the plane}\label{8855}

Now we go back to the regularity issue of the $s$-minimal sets.
Since this topic seems to be very difficult to deal with in the general case,
we start with the case of low dimension~$n=2$.
For this, first 
we point out 
that 
\begin{equation}\label{K not}
{\mbox{the cone~${\mathcal{K}}$ in~\eqref{cone.1},
that is the ``black cone'' in Figure~3,
is not $s$-minimal.}}\end{equation}
The proof given here is due to an original idea of L. Caffarelli.
Suppose, by contradiction, that ${\mathcal{K}}$
is $s$-minimal. Then 
consider the 
set~${\mathcal{K}}'$
in Figure~5 that is obtained from ${\mathcal{K}}$
by adding another little square adjacent to the origin.
\medskip

\begin{center} \label{cone2}
\includegraphics[width=4in]{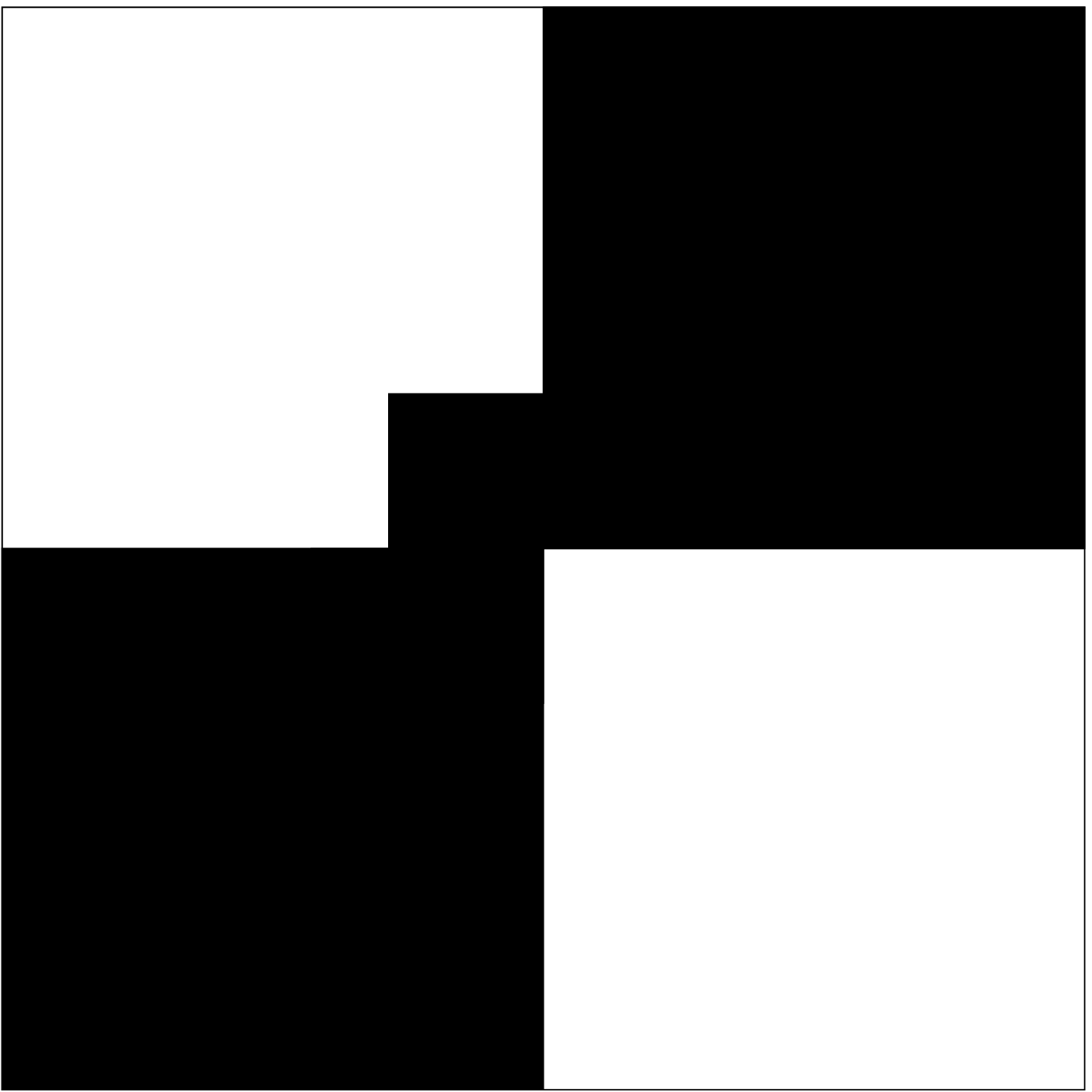}
\end{center}
\nopagebreak\centerline{\footnotesize\em
Figure~5.
The cone ${\mathcal{K}}'$.}
\bigskip

Then, the $s$-perimeter of~${\mathcal{K}}'$ (in a domain
large enough to contain the additional
little square)
is equal to that of~${\mathcal{K}}$.

To check this, just compare Figures~6
and~7: in~${\mathcal{K}}$ the additional little square
is ``white'' and therefore it interacts with the ``black quadrants''
$A$ and $B$ and with the ``black rectangles'' $C$ and $D$
in Figure~6,
while in~${\mathcal{K}}'$ the additional little square
is ``black'' and therefore it interacts with the ``white quadrants''
$A'$ and $B'$ and with the ``white rectangles'' $C'$ and $D'$.
Since the interactions with~$A\cup B$ (resp.,
$C\cup D$) are equal to the ones with~$A'\cup 
B'$ (resp., $C'\cup D'$), we have that
the $s$-perimeter of~${\mathcal{K}}'$ 
is equal to that of~${\mathcal{K}}$
(notice that -- due to the finite space at our disposal --
Figures~6 and~7 only represent a ``bounded portion'' of~$\R^2$,
and the sets~$A$, $B$, $C$, $D$, $A'$, $B'$, $C'$ and~$D'$
are actually all unbounded).
\medskip

\begin{center} \label{cone3}
\includegraphics[width=4in]{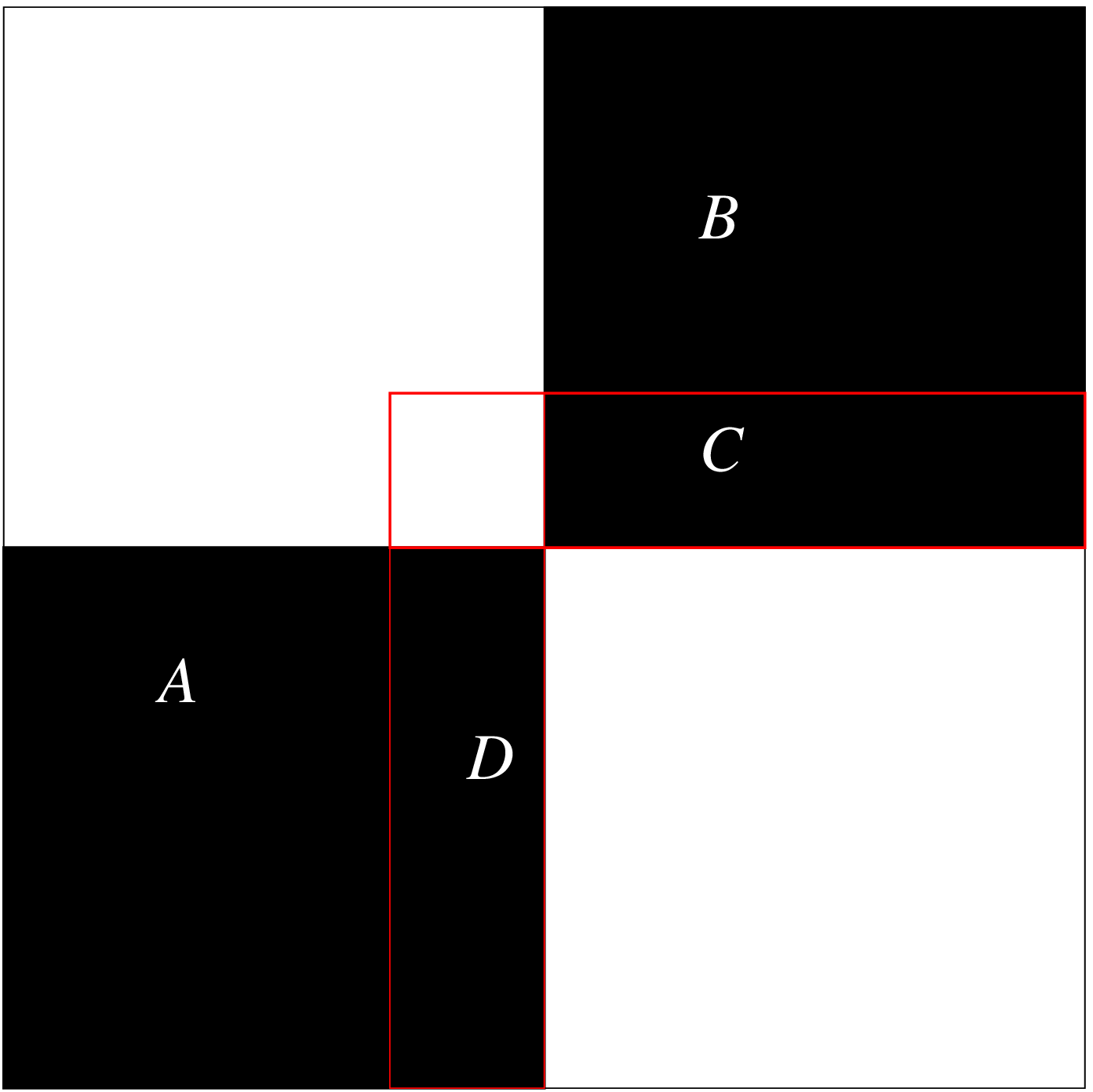}
\end{center}
\nopagebreak\centerline{\footnotesize\em
Figure~6. The sets~$A$, $B$, $C$ and~$D$
that interact with the little white square in ${\mathcal{K}}$.}
\bigskip

\begin{center} \label{cone4}
\includegraphics[width=4in]{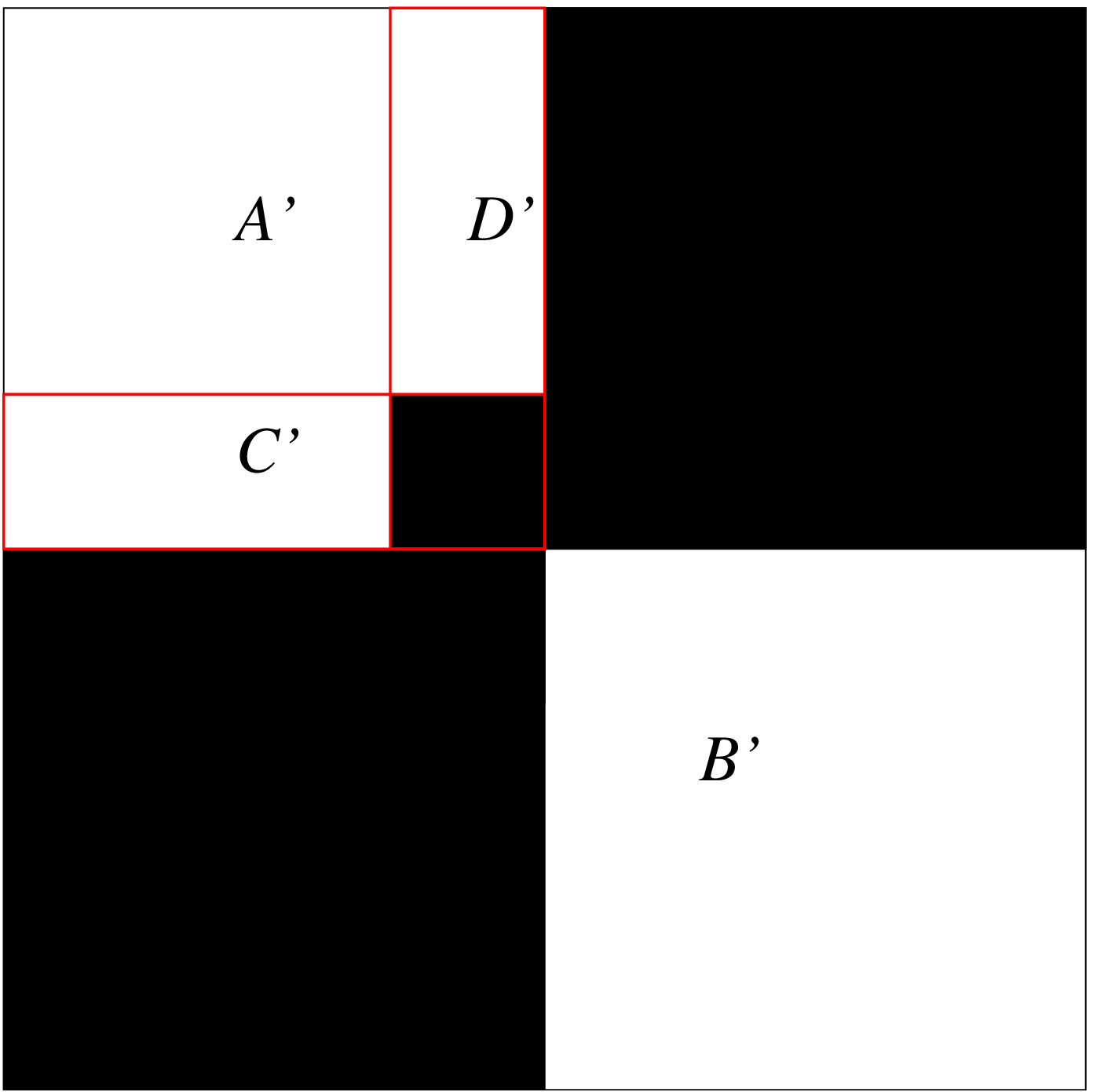}
\end{center}
\nopagebreak\centerline{\footnotesize\em
Figure~7.
The sets~$A'$, $B'$, $C'$ and~$D'$ 
that interact with the little black square in ${\mathcal{K}}'$.}
\bigskip

As a consequence,~${\mathcal{K}}'$ is $s$-minimal
(since we supposed that so is~${\mathcal{K}}$),
and therefore~${\mathcal{K}}'$ satisfies the Euler-Lagrange 
equation in~\eqref{EL} at the origin.
But this cannot be, since the ``black region''
contributes more than the ``white one'', namely
$$ \myint_{\R^n}\frac{\chi_{{ \mathcal{K}}'}(y)-
\chi_{\R^n\setminus{ \mathcal{K}}'}(y)
}{|y|^{n+2s}}\,dy >0.$$
This contradiction shows that the cone~${ \mathcal{K}}$
is not $s$-minimal, proving~\eqref{K not}.

Moreover, recalling~\eqref{K el}, we have that~${ \mathcal{K}}$
is an example of a set that satisfies~\eqref{EL} but that is not 
$s$-minimal: thus the Euler-Lagrange equation in~\eqref{EL}
is implied by, but it is not equivalent to, $s$-minimality.

It would be interesting to construct examples (if they exist)
of smooth sets that satisfy 
the Euler-Lagrange equation in~\eqref{EL} without being
$s$-minimal.

{F}rom~\eqref{K not} one may conjecture that $s$-minimal sets
are smooth in dimension~$n=2$ (I mean, if any singularity occurs,
one can prove that one can reduce to
a cone, and so one should suspect that 
the ``worst'' cone is the ninety degree one in Figure~3).
Unfortunately it is not easy to extend the above geometric argument
to the general case (for instance, the singular cone could be
made of many sectors, and these sectors could differ
one from the other, see Figure~10). To get around this difficulty, 
in~\cite{SV3}
the regularity of $s$-minimal sets in dimension~$2$ is proved
using an analytic argument based on domain perturbations. The
result obtained\footnote{In the regularity results
such as Theorem~\ref{2}, we are implicitly ruling out the trivial cases
in which either~$E$ or its complement is empty.
Also, for the sake of precision, we point out that
in~\cite{SV3} the regularity obtained is only
of~$C^{1,\alpha}$-type:
the bootstrap improvement towards~$C^\infty$-regularity
is not trivial and it is contained \label{d77duuquqq}
in~\cite{barrios}.}
is the following:

\begin{theorem}[\cite{SV3}]\label{2}
Let~$n=2$. If~$R>r>0$ and~$E$ is an $s$-minimal set in~$B_R$, 
then~$(\partial E) 
\cap B_r$ is a $C^\infty$-curve.
If $E$ is an $s$-minimal set in $B_\rho$ for every~$\rho>0$, 
then~$\partial E$ is 
a straight line.
\end{theorem}

As a byproduct of Theorem~\ref{2} and of
a dimensional reduction in~\cite{CRS},
one also improves Theorem~\ref{s989}, obtaining that the
singular set $\Sigma$ in~$\R^n$
has finite~$(n-3)$-Hausdorff dimension (instead of~$n-2$:
and we do not know whether or not this is optimal, see Theorem~\ref{n} below).

The last claim in Theorem~\ref{2} somehow states that fractional
geodesics in the plane are straight lines, as happens in the 
classical case. The proof of Theorem~\ref{2} is based on domain
perturbation. The idea of the proof may be sketched by
thinking about classical geodesics in the plane.
\medskip

\begin{center} \label{geo1}
\includegraphics[width=4in]{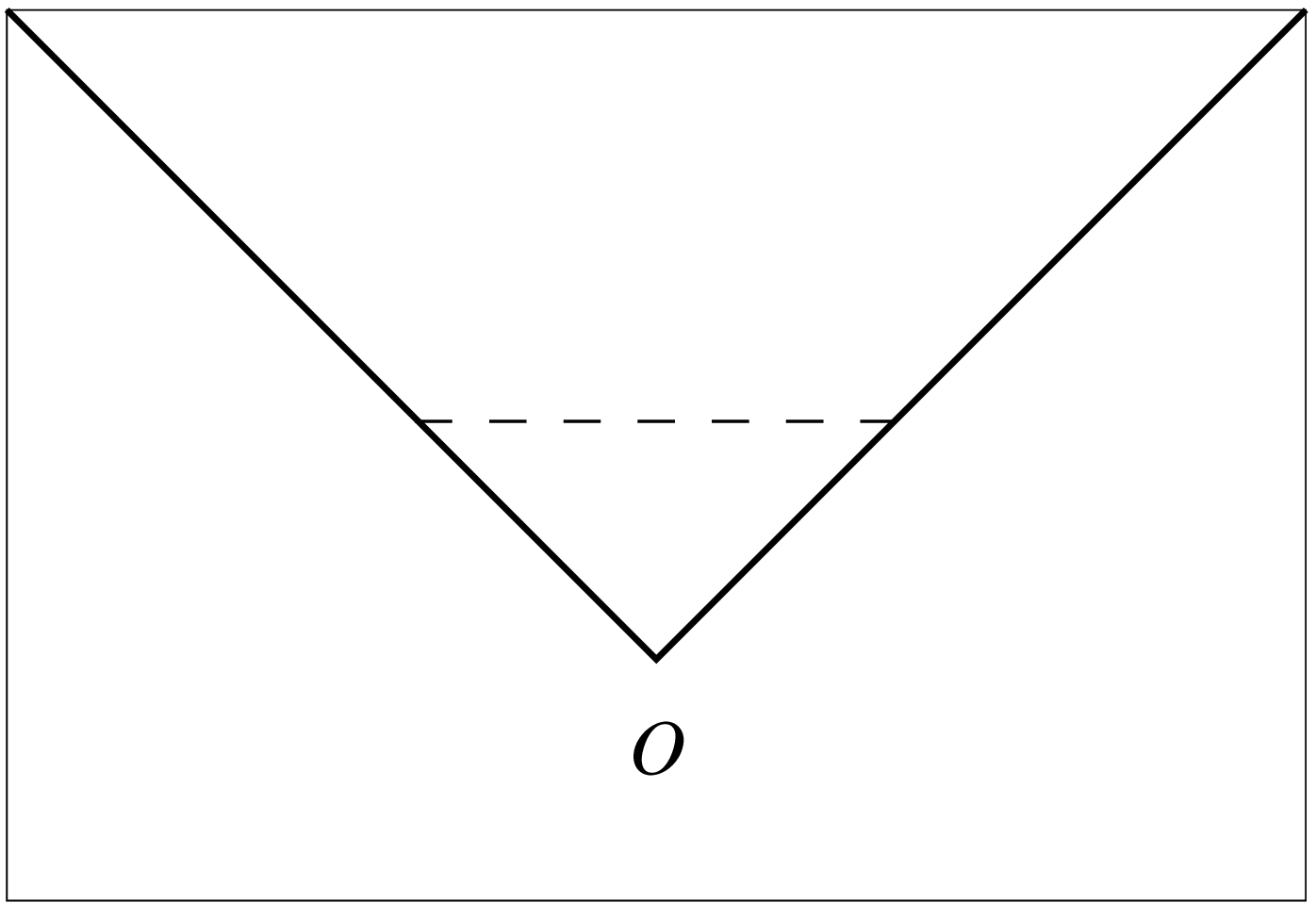}
\end{center}
\nopagebreak\centerline{\footnotesize\em
Figure~8.
An edge in the plane is not a geodesic: the classical proof}
\bigskip

The classical proof to show that an edge
is not a geodesic consists in cutting the angle in~$O$
and shortening the length by the dashed segment
as shown in Figure~8. This type of proof
is difficult to transpose into a fractional framework, since
the new object is not a smooth deformation of the original one.
But there is a modification of this argument that shortens the
length by taking a domain perturbation near~$O$ of the edge
plus a suitable gluing at infinity. This alternative argument
is depicted in Figure~9.
\medskip

\begin{center} \label{geo2}
\includegraphics[width=4in]{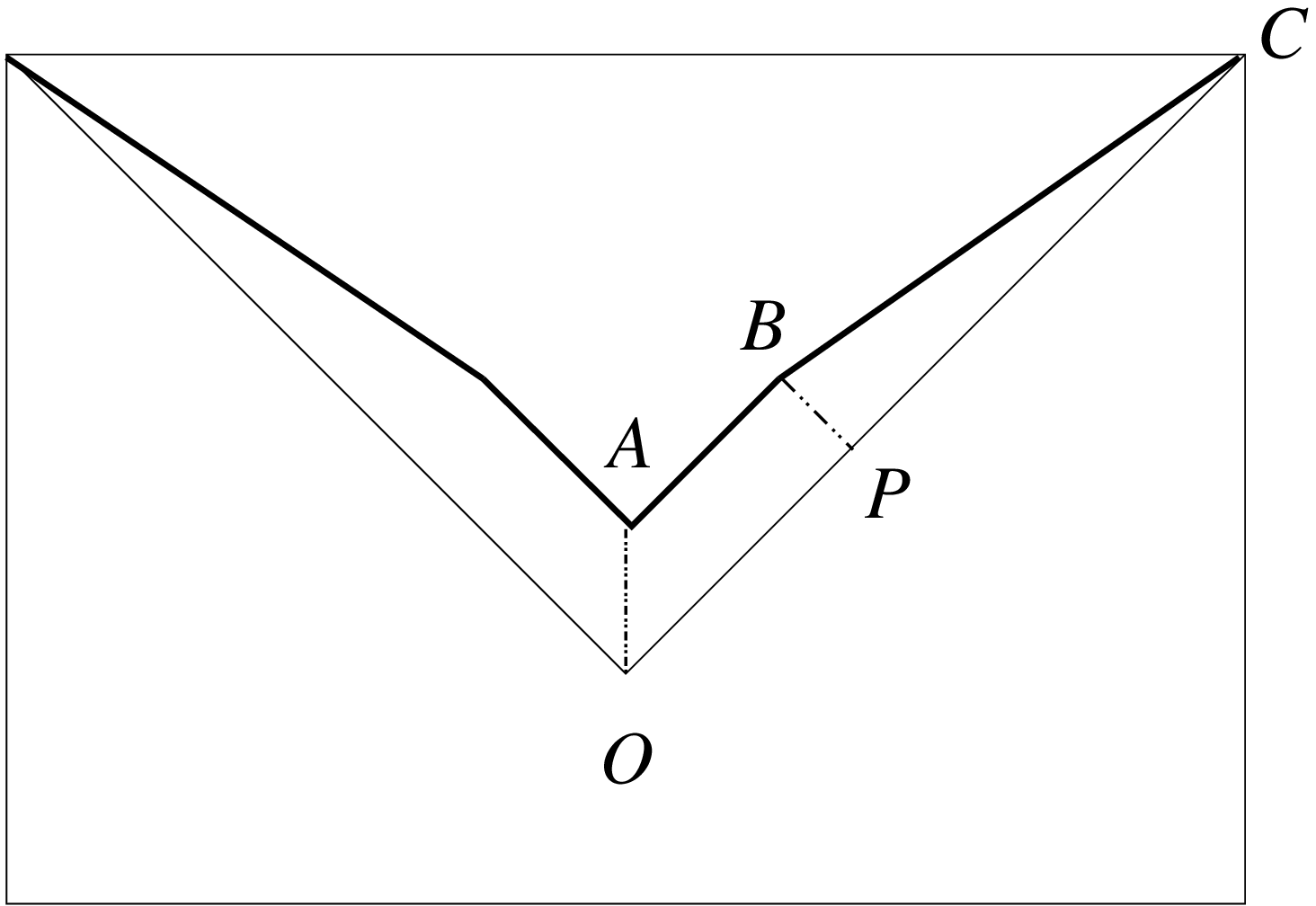}
\end{center}
\nopagebreak\centerline{\footnotesize\em
Figure~9.
An edge in the plane is not a geodesic: the domain variation
proof.}
\bigskip

The argument goes like this (we follow the right side of
the picture, the left one being symmetric). We translate
the vertex slightly upwards, say, in such a way that~$AO$
has length~$1$. Then, the length of~$AB$ is strictly shorter than the
one of~$OP$, say
\begin{equation}\label{d}
\overline{AB}=
\overline{OP}-\delta\end{equation}
for some~$\delta>0$. This is 
not
a contradiction yet, since $AB$ is not a compact modification of~$OP$,
so, for our purposes, we need to glue~$AB$ with~$OP$. For this, we 
take a suitably large $R>0$, and we join~$B$ to the point~$C$, which
is chosen in the half line from~$O$ to~$P$ in such a way that~$PC$
has length~$R$. Then, by Pythagoras' Theorem,
\begin{equation}\label{74}
\overline{BC}=\sqrt{\overline{BP}^2+\overline{PC}^2}\le
\sqrt{\overline{AO}^2+\overline{PC}^2}=\sqrt{1+R^2}\le
R+\frac{C}{R}.\end{equation}
Now, if~$R$ is chosen large enough, we obtain that
the polygonal chain~$ABC$ is shorter than the segment~$OC$, namely:
\begin{equation}\label{over}
\overline{OC}-\overline{ABC}=
\overline{OP}+\overline{PC}-\big(\overline{AB}
+\overline{BC}\big)\ge \delta+R-\Big(R+\frac{C}{R}\Big)
=\delta-\frac{C}{R}>0.\end{equation}
The argument on the left in Figure~9
is the same, and so we have shown that the edge has a longer length
than the polygonal chain on the top of it.

This alternative argument proving that the edges in the plane
are not geodesic may be extended to the fractional case. That is, one considers
a minimal cone~$E\subset \R^2$ (different than a hyperplane)
and constructs a set $\tilde E$ as a translation 
of $E$ in $B_{R/2}$ which coincides with $E$ outside $B_R$. 
Then we use Proposition~\ref{73} to compute
the difference between the energies of the extensions of $\tilde E$ 
and $E$: we obtain that this difference is bounded by~$C/R^{2s}$
(notice that formally the limit case~$s=1/2$ goes back to
the term~$C/R$ in~\eqref{74}).
Technically, this estimate is achieved by considering a smooth
perturbation~$\phi\in C_0^\infty(\R^{3})$,
with~$\phi(X)=1$ if~$|X|\le 1/2$ and~$\phi(X)=0$ if~$|X|\ge 3/4$, and
considering the diffeomorphism
\begin{equation} \label{34}
\R^3\ni X\mapsto Y(X) :=X+ \Big(\phi(X/R),0,0\Big).\end{equation}
The inverse diffeomorphism is denoted, with a slight but common
abuse of notation, $\R^3\ni Y\mapsto X(Y)=X$.
Given~$u$ to be the extension of~$\chi_E$, according 
to Proposition~\ref{73}, one
defines~$u^+_R(Y):=u(X)$.
Similarly, by switching~$\phi$
with~$-\phi$ in~\eqref{34}, we can define~$u^-_R(Y)$.
Of course, the derivatives of $u^\pm$ may be computed
from the ones of~$u$ via the Chain Rule: in this way,
we can compute~${\mathcal{E}}_{B_R}(u^\pm)$ and compare it 
with~${\mathcal{E}}_{B_R}(u)$: one obtains
$$ {\mathcal{E}}_{B_R}(u^+_R)+{\mathcal{E}}_{B_R}(u^-_R)-2
{\mathcal{E}}_{B_R}(u) \le 
\frac{C}{R^{2s}} .$$
Then, since $u$ is a minimiser, ${\mathcal{E}}_{B_R}(u)\le
{\mathcal{E}}_{B_R}(u^-_R)$ and so we obtain
\begin{equation} \label{X}
{\mathcal{E}}_{B_R}(u^+_R)-{\mathcal{E}}_{B_R}(u)\le
\frac{C}{R^{2s}} .\end{equation}
Now we look at the cone~$E$ in~$\R^2$: up
to a rotation, we may suppose that a sector of~$E$ has an angle less
than~$\pi$ and
is bisected by~$e_2$. Thus,
there exist~$M\ge1$ and~$p\in B_M$, on the $e_2$-axis, such that $p$
lies in the interior of $E$,
and~$p+ e_1$ and $p-e_1$ lie in the exterior of $E$,
and we let~$P:=(p,0)\in\R^3$ (see Figure~10
where~$q:=p+e_1$).
\medskip

\begin{center} 
\includegraphics[width=4in]{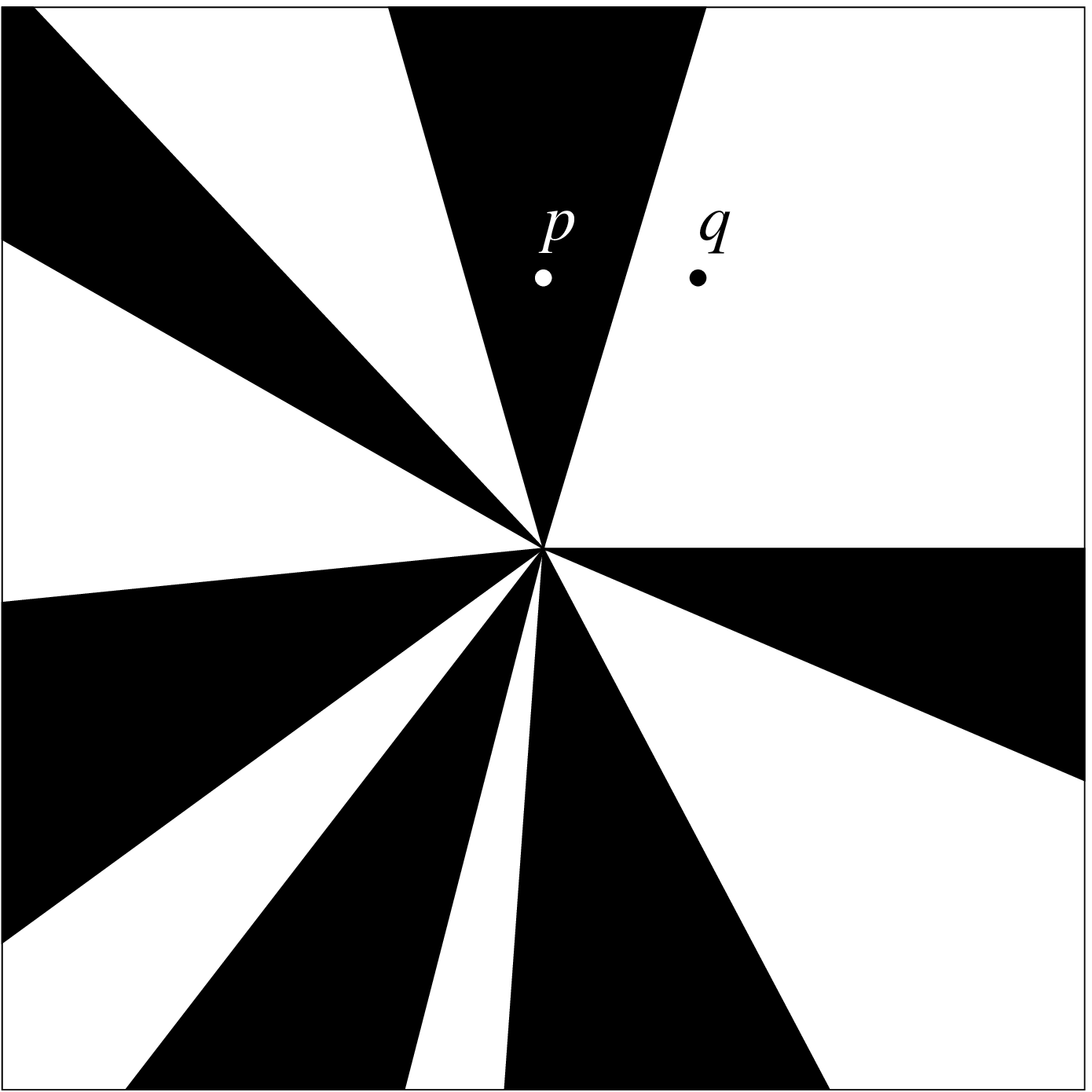}
\end{center}
\nopagebreak\centerline{\footnotesize\em
Figure~10.
The singular cone $E$, with~$p\in B_M$ and~$q:=p+e_1$.}
\bigskip

Taking~$R$ much larger than~$M$ we see that
$u^+_R(Y)=u(Y-e_1)$
if~$|Y|\le 2M$, and~$u^+_R$ coincides with~$u$ if~$|Y|\ge R$.
We define 
\begin{equation}\label{1.14a}
v_R(X):=\min \{ u(X), \,u^+_R(X)\} \quad  {\mbox{ and }} \quad 
w_R(X):=\max \{ u(X), \,u^+_R(X)\}.\end{equation}
By construction,
$u^+_R<w_R=u$
in a neighbourhood of~$P$, and
$u<w_R=u^+_R$
in a neighbourhood of~$P+e_1$, that is~$u$ and~$u^+_R$ cross
each other inside the ball of radius~$2M$. This and the maximum principle
imply that~$w_R$ (as well as~$v_R$) cannot be a minimiser with
respect to compact perturbations in the ball of radius~${2M}$: that is, there 
exists~$\delta>0$ 
and a modification~$u_*$ of~$w_R$ inside~$B_{2M}$
such that
\begin{equation}\label{92}
{\mathcal{E}}_{B_{2M}}(u_*)\le {\mathcal{E}}_{B_{2M}}(w_R)-\delta.
\end{equation}
Notice that this~$\delta>0$ is independent of~$R$
(since~$w_R$ restricted to the ball of radius~${2M}$
is simply the
maximum between~$u$ and its translation and so it is
is independent of~$R$). The role played by~$\delta$ here
is indeed analogous to the one of~\eqref{d}.
Since~$u_*=w_R$ outside~$B_{2M}$ we have that
$$ {\mathcal{E}}_{B_R\setminus
B_{2M}}(u_*)={\mathcal{E}}_{B_R\setminus B_{2M}}(w_R),$$
and so~\eqref{92} becomes
\begin{equation}\label{92.1} 
{\mathcal{E}}_{B_{R}}(u_*)\le {\mathcal{E}}_{B_{R}}(w_R)-\delta.
\end{equation}
The advantage of working with~\eqref{92.1}
rather than~\eqref{92} is that the energy domain is now the ball
of radius~$R$ (that is the domain that supports the perturbation),
but~$\delta$ is independent of~$R$.
The minimality of~$u$ also gives that
\begin{equation}\label{92.2} 
{\mathcal{E}}_{B_{R}}(u)\le
{\mathcal{E}}_{B_{R}}(v_R) .
\end{equation}
Now, in light of~\eqref{1.14a}, we point out the integral identity
\begin{equation}\label{92.3} 
{\mathcal{E}}_{B_{R}}(v_R)+{\mathcal{E}}_{B_{R}}(w_R)=
{\mathcal{E}}_{B_{R}}(u)+{\mathcal{E}}_{B_{R}}(u^+_R).\end{equation}
All in all, we have that
\begin{equation}\label{over2}
\begin{split}
& {\mathcal{E}}_{B_{R}}(u)-
{\mathcal{E}}_{B_{R}}(u_*)\\
{\mbox{by \eqref{92.1} }}&\qquad\ge\,
{\mathcal{E}}_{B_{R}}(u)-
{\mathcal{E}}_{B_{R}}(w_R)+\delta
\\
{\mbox{by \eqref{92.3} }}&\qquad=\,
{\mathcal{E}}_{B_{R}}(v_R)-
{\mathcal{E}}_{B_{R}}(u^+_R)+\delta
\\
{\mbox{by \eqref{92.2} }}&\qquad\ge\,
{\mathcal{E}}_{B_{R}}(u)-
{\mathcal{E}}_{B_{R}}(u^+_R)+\delta
\\
{\mbox{by \eqref{X} }}&\qquad\ge\,
\delta-\frac{C}{R^{2s}},
\end{split}\end{equation}
which is strictly positive for $R$ large enough.
This is in contradiction with the minimality of~$u$ and so
it completes the proof of Theorem~\ref{2}.

Notice that~\eqref{over2} plays the role of~\eqref{over}
in this case. Furthermore, the technique used to prove Theorem~\ref{2}
seems to work for a wide class of variational problems:
see e.g.~\cite{SV gen}, where these ideas are exploited to prove
monotonicity and symmetry results for minimisers and stable solutions
of a very general class of functionals.

\subsection{Regularity of $s$-minimal sets when $s$ is close 
to~$1/2$}\label{88811}

Having settled the regularity of~$s$-minimal sets in
the plane in Theorem~\ref{2},
we discuss now the possible regularity in a higher dimension.
As far as we know, this problem is
open up to now. Though no example of a singular set in any
dimension and for any~$s\in(0,1/2)$ is known, the only regularity
result available at the moment seems to be the following one,
which recovers the classical minimal surface regularity theory
when~$s$ is sufficiently close\footnote{As already pointed out
in the footnote on page~\pageref{d77duuquqq}, in~\cite{CV2}
the regularity theory is of~$C^{1,\alpha}$-type: for the bootstrap
to~$C^\infty$-regularity see~\cite{barrios}.}
to~$1/2$:

\begin{theorem}[\cite{CV2}]
\label{n}
For any~$n\in\N$ there exists~$\epsilon_n\in (0,1/2]$ such
that if~$s\in ((1/2)-\epsilon_n,\,1/2)$ then $s$-minimal sets
are ``as regular as the classical minimal surfaces in dimension~$n$'',
namely:
\begin{itemize}
\item If $n\le 7$ and~$s\in ((1/2)-\epsilon_n,\,1/2)$, then
any $s$-minimal set is locally a $C^\infty$-surface.
\item If $n=8$ and~$s\in ((1/2)-\epsilon_8,\,1/2)$, then
any $s$-minimal set is locally a $C^\infty$-surface
except, at most, at countably many isolated points.
\item If $n\ge8$ and~$s\in ((1/2)-\epsilon_n,\,1/2)$, then
any $s$-minimal set is locally a $C^\infty$-surface
outside a closed set~$\Sigma\subset\partial E$ with
finite~$(n-8)$-Hausdorff dimension.
\end{itemize}
\end{theorem}

Of course, in the notation of Theorem~\ref{n}, $\Sigma$ could well be
empty. The finite~$(n-8)$-Hausdorff dimension statement
in Theorem~\ref{n} improves (when~$s$ in the ``right range'')
the previous
ones mentioned in
Theorem~\ref{s989} and right below Theorem~\ref{2}.
Unfortunately the proof of Theorem~\ref{n} uses a compactness argument,
therefore nothing is known on~$\epsilon_n$
(except that it is a positive, universal quantity, depending only on~$n$,
but no explicit bound is available). Of course, from Theorem~\ref{2}
we know that~$\epsilon_2=1/2$ is fine for the regularity theory when~$n=2$
(but this really comes from~\cite{SV3}
and it cannot be proved with the argument in~\cite{CV2}).
Of course, any explicit bound on~$\epsilon_n$ would be welcome
to make Theorem~\ref{n} applicable in concrete cases.

Notice also that, in view of Theorem~\ref{ts 0},
we know that for~$s$ close to~$0$ the $s$-minimal sets
seem related to the minimisers of the Lebesgue measure,
for which no regularity result is possible
(a set can have a small Lebesgue measure
and possess a very wild boundary). Therefore,
the regularity of $s$-minimal sets when $s$ is close to~$0$
(if it holds true) is conceivably more difficult
than in the case
in which $s$ is close to~$1/2$.

\section{The fractional Allen-Cahn equation}

Classical minimal surfaces naturally arise in phase transition models.
Similarly $s$-minimisers of the functional in~\eqref{P}
arise in non-local phase transition models, in which
the classical diffusion term is replaced by the fractional Laplacian.
To see this, let us briefly review the relation between phase transitions
and minimal surfaces in the standard case. We take~$W\in C^2(\R)$ to be a 
``double-well
potential'', say, for concreteness,
$$ W(t):=\frac{(1-t^2)^2}{4}.$$
Then, the classical Allen-Cahn (or scalar Ginzburg-Landau)
phase coexistence model investigates the functional
\begin{equation}\label{0F}
{\mathcal{F}}(u,U):=\myint_U \frac{|\nabla 
u(x)|^2}2+W(u(x))\,dx.\end{equation}
The minimisers of this
functional satisfy the
Allen-Cahn equation
\begin{equation}\label{AC}
-\Delta u=u-u^3 \,{\mbox{ in }}\, U,
\end{equation}
and they
have a strong tendency
to stay close to~$\pm 1$ (which are the ``pure
phases'' of the model) since these values kill the 
potential energy, while the gradient term forces
the transition to occur with the least possible surface tension.
These heuristic considerations can be made rigorous by introducing
a parameter~$\epsilon$ and by considering the rescaled 
minimiser
\begin{equation}\label{ue}
u_\epsilon(x):= u(x/\epsilon).\end{equation}
Scaling~$u$ to $u_\epsilon$ in~\eqref{0F}
(and normalising by a multiplicative factor of order~$\epsilon^{n-1}$ which
does not change the notion of minimisers), one is lead to study the
functional
\begin{equation}\label{0Feps}
{\mathcal{F}}_\epsilon(u,U):=\myint_U \frac{\epsilon |\nabla
u(x)|^2}2+\frac{1}{\epsilon} W(u(x))\,dx.\end{equation}
Then, the following classical result holds true:

\begin{theorem}[\cite{mod,cord}]\label{sd9fffffff}
$\,$
\begin{itemize}
\item As~$\epsilon\searrow0$, ${\mathcal{F}}_\epsilon$ 
$\Gamma$-converges
to the classical perimeter functional, i.e., for any set~$E$
of locally finite perimeter, the following holds:
\begin{itemize}
\item For any $u_\epsilon\in L^1_{\rm loc}(\R^n,[-1,1])$
converging to~$\chi_E-\chi_{\R^n\setminus E}$ in~$L^1_{\rm loc}(\R^n)$, 
we have that
$$ \Per (E,U)\le \liminf_{\epsilon\searrow0}{\mathcal{F}}_\epsilon
(u_\epsilon,U);$$
\item There exists~$u_\epsilon\in L^1_{\rm loc}(\R^n,[-1,1])$ that
converges to~$\chi_E-\chi_{\R^n\setminus E}$ in~$L^1_{\rm loc}(\R^n)$ 
such that
$$ \Per (E,U)\ge \limsup_{\epsilon\searrow0}{\mathcal{F}}_\epsilon
(u_\epsilon,U);$$
\end{itemize}
\item The following compactness properties holds:
if $u_\epsilon\in L^1_{\rm loc}(\R^n,[-1,1])$ and
$$ \sup_{\epsilon\in(0,1)}
{\mathcal{F}}_\epsilon(u_\epsilon,U)<+\infty,$$
then there exists~$E$ and a convergent subsequence such 
that~$u_\epsilon$
converges to~$\chi_E-\chi_{\R^n\setminus E}$ in~$L^1_{\rm loc}(\R^n)$.
\item Fixed~$R>r>0$, $\vartheta_1$, $\vartheta_2\in(-1,1)$, 
if~$u_\epsilon$ 
minimises ${\mathcal{F}}_\epsilon$
in~$B_R$ (i.e. 
if
${\mathcal{F}}(u,B_R)\le {\mathcal{F}}_{B_R}(u+\varphi)$
for any~$\varphi\in C^\infty_0(B_R)$), and $u_\epsilon(0)>\vartheta_1$
then
$$ {\mathcal{L}}^n\Big( B_R\cap \{u_\epsilon>\vartheta_2\}\Big)
\ge cR^n,$$
provided that~$\epsilon< c(\vartheta_1,\vartheta_2)R$.
Also, $\{u_\epsilon>\vartheta_2\}$ 
approaches~$E$ uniformly in~$B_r$, and~$E$ minimises the perimeter
in~$B_r$ with respect to its boundary data.\end{itemize}
\end{theorem}

The aim of the following pages is to present what happens to these
results in the fractional framework and to discuss some possible
consequences. For this, we first introduce a domain notation
by setting
\begin{eqnarray*}
Q_U &:=&\big(U\times U\big)\cup \big((\R^n\setminus 
U)\times U\big)
\cup\big(U\times (\R^n\setminus U)\big)\\
&=&\R^{2n}\setminus \big(
(\R^n\setminus U)\times(\R^n\setminus U)
\big).\end{eqnarray*}
The set~$Q_U$ will represent the natural domain of
a non-local interaction between particles in~$\R^n$:
namely this interaction is represented by an integral 
over~$\R^n\times\R^n$, but we will remove from this integral the
contribution given by two particles both lying in the complement of~$U$,
since this will be considered fixed by the datum
(this is the same type of renormalisation
procedure performed in~\eqref{P}).
More concretely, for any~$s\in(0,1)$ we consider the functional
\begin{equation}\label{0G}
{\mathcal{G}}(u;U):=
\myiint_{Q_U}\frac{|u(x)-u(y)|^2}{2\,|x-y|^{n+2s}}\,dx\,dy+
\myint_U W(u(x))\,dx.\end{equation}
Notice that the functional~${\mathcal{G}}$ differs from
the functional~${\mathcal{F}}$ in~\eqref{0F}
since the gradient part (i.e., the $H^1$-seminorm of~$u$ in~$U$)
is replaced here by a double integral of Gagliardo type,
which tries to mimic a long-range particle interaction energy.
The Euler-Lagrange equation associated with~${\mathcal{G}}$ is
\begin{equation}\label{AC2}
(-\Delta)^s u=u-u^3 \,{\mbox{ in }}\, U,\end{equation}
which may be seen as a fractional variant of the classical
Allen-Cahn equation in~\eqref{AC} (see, e.g.,
\cite{Im} for applications to fractional
mean curvature flows). To obtain
a functional on which a $\Gamma$-convergence problem is well-posed,
we proceed as before, taking~$u_\epsilon$ as in~\eqref{ue},
and scaling~$u$ to $u_\epsilon$ in~\eqref{0G}:
unlike the classical case, here it is also necessary
to normalise by a multiplicative factor that depends on~$s$,
namely~$\epsilon^{n-2s}$ if~$s\in(0,1/2)$, 
$\epsilon^{n-1}$ if~$s\in(1/2,1)$
and~$\epsilon^{n-1}\log(1/\epsilon)$ when~$s=1/2$
(often, in fractional problems, a logarithmic correction at~$s=1/2$
is necessary to match the case~$s\in(0,1/2)$ with the 
case~$s\in(1/2,1)$). This procedure leads to the following functional
$${\mathcal{G}}_\epsilon(u;U):=\left\{
\begin{matrix}
\displaystyle\myiint_{Q_U}\frac{|u(x)-u(y)|^2}{2\,|x-y|^{n+2s}}\,dx\,dy+
\displaystyle\frac{1}{\epsilon^{2s}}
\displaystyle\myint_U W(u(x))\,dx, & {\mbox{ if 
}}s\in(0,1/2),\\
\\
\displaystyle\log(1/\epsilon) 
\displaystyle\myiint_{Q_U}\frac{|u(x)-u(y)|^2}{2\,|x-y|^{n+2s}}\,dx\,dy+
\displaystyle\frac{\log(1/\epsilon)}{\epsilon}
\displaystyle\myint_U W(u(x))\,dx, & {\mbox{ if }}s=1/2,\\
\\
\epsilon^{2s-1}
\displaystyle\myiint_{Q_U}\frac{|u(x)-u(y)|^2}{2\,|x-y|^{n+2s}}\,dx\,dy+
\displaystyle\frac{1}{\epsilon}
\displaystyle\myint_U W(u(x))\,dx, & {\mbox{ if }}s\in(1/2,1).
\end{matrix}\right. $$

\subsection{$\Gamma$-convergence and density estimates}

In this non-local setting for phase transitions,
one may recover the $\Gamma$-convergence and density estimates
results of Theorem~\ref{sd9fffffff}, though they now
provide an alternative selection of local
and non-local limit interfaces, according to the cases~$s\in(0,1/2)$
and~$s\in[1/2,1)$. Namely, the limit problem reduces to
the $s$-perimeter functional\footnote{The careful reader
will have noticed that the $s$-perimeter functional is defined only
for~$s\in(0,1/2)$ while the fractional Allen-Cahn
equation for any~$s\in(0,1)$.}
when~$s\in(0,1/2)$
and to the classical perimeter functional when~$s\in[1/2,1)$.
This is somehow consistent with the fact that the $s$-perimeter
reduces to the classical perimeter as~$s\nearrow1/2$
(recall Theorem~\ref{888990}). Also, from the point of view
of the applications, it suggests that the limit interfaces
of the non-local Allen-Cahn phase transition may be either
local or non-local according to whether the fractional parameter~$s$
is above or below the critical threshold~$1/2$ (that is, when~$s\ge1/2$
the nonlocal effect is lost by the limit interface). In further
detail, the result obtained reads as follows:\footnote{Some preliminary
work needed for the proof of Theorem~\ref{d88fdghjww} and a careful
analysis of the
one-dimensional case was also performed in~\cite{palsa}.}

\begin{theorem}[\cite{SV1, SV Sob, SV2}]\label{d88fdghjww}
$\,$
\begin{itemize}
\item As~$\epsilon\searrow0$, ${\mathcal{G}}_\epsilon$ 
$\Gamma$-converges
to the $s$-perimeter functional when~$s\in(0,1/2)$ and
the classical perimeter functional when~$s\in[1/2,1)$
\item The following compactness property holds:
if $u_\epsilon\in L^1_{\rm loc}(\R^n,[-1,1])$ and
$$ \sup_{\epsilon\in(0,1)}
{\mathcal{G}}_\epsilon(u_\epsilon,U)<+\infty,$$
then there exists~$E$ and a convergent subsequence such 
that~$u_\epsilon$
converges to~$\chi_E-\chi_{\R^n\setminus E}$ in~$L^1_{\rm loc}(\R^n)$.
\item Fixed~$R>r>0$, $\vartheta_1$, $\vartheta_2\in(-1,1)$, 
if~$u_\epsilon$ 
minimises ${\mathcal{G}}_\epsilon$
in~$B_R$ 
and $u_\epsilon(0)>\vartheta_1$
then
$$ {\mathcal{L}}^n\Big( B_R\cap \{u_\epsilon>\vartheta_2\}\Big)
\ge cR^n,$$
provided that~$\epsilon< c(\vartheta_1,\vartheta_2)R$.
Also, $\{u_\epsilon>\vartheta_2\}$ 
approaches~$E$ uniformly in~$B_r$, and~$E$ minimises either the 
$s$-perimeter or the classical perimeter
in~$B_r$ with respect to its boundary data
(depending on whether~$s\in(0,1/2)$ or~$s\in[1/2,1)$).\end{itemize}
\end{theorem}

We try to translate the statement of Theorem~\ref{d88fdghjww}
into an evocative picture that involves the 
parameters~$\epsilon$
and~$s\in(0,1)$. Namely, in Figure~11, $s$ ranges horizontally
and~$\epsilon$ vertically; on the top of the picture
(corresponding to the case~$\epsilon=1$)
we have a phase transition function~$u_\epsilon$ whose level
sets as $\epsilon\searrow0$ approach some~$\partial E$,
which is drawn in the bottom of the picture
(which corresponds to the case~$\epsilon=0$).
When~$s\in(0,1/2)$ this~$\partial E$ is an~$s$-minimal set,
while for~$s\in[1/2,1)$ it is a classical minimal surface
(since ``classical minimal surfaces are nice''
and ``$s$-minimal sets might be somewhat wild''
the picture is trying to distinguish between them
by showing either smooth surfaces or singular cones).
\medskip

\begin{center} 
\includegraphics[width=4in]{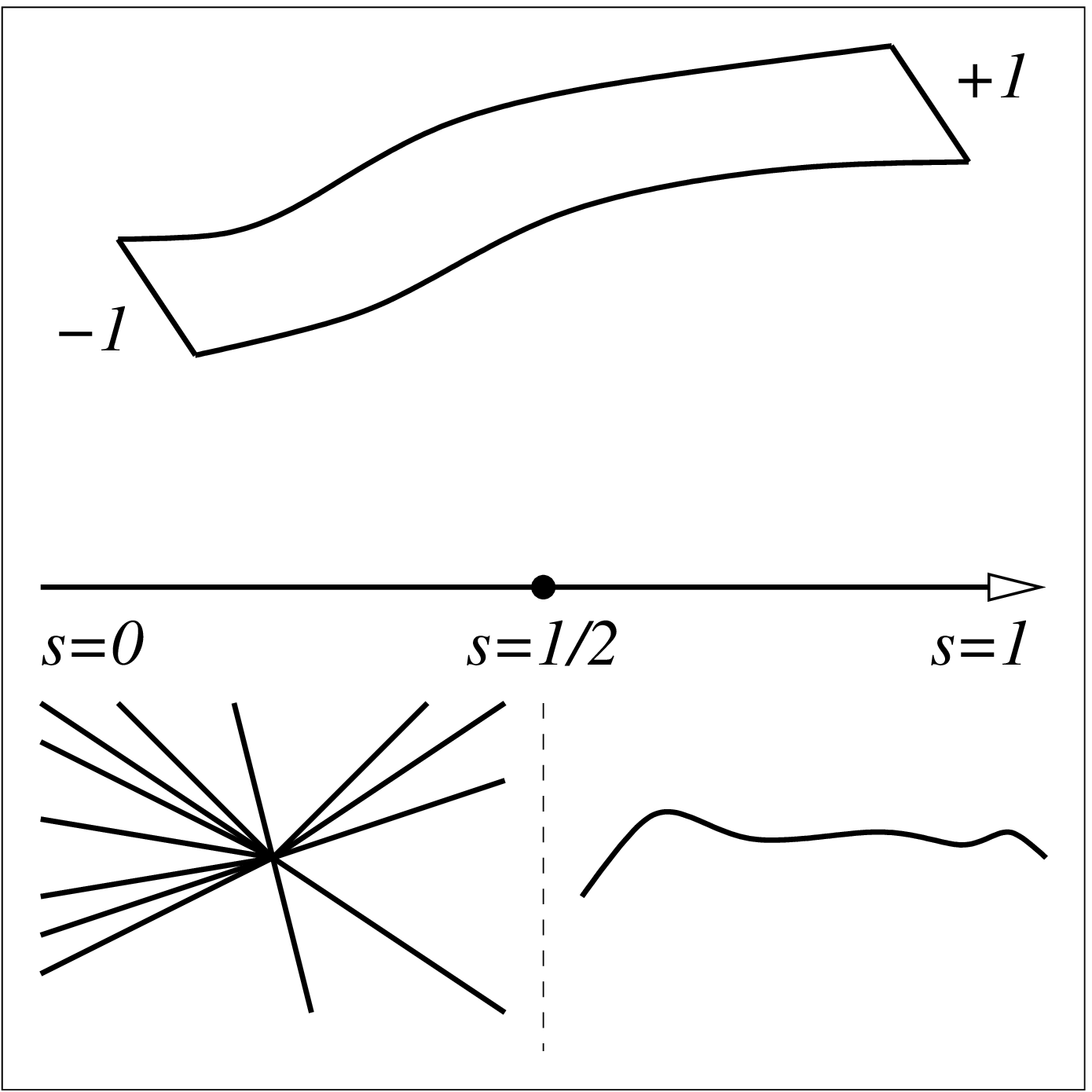}
\end{center}
\nopagebreak\centerline{\footnotesize\em
Figure~11. $\Gamma$-convergence for fractional
phase transitions.}
\bigskip

As a matter of fact, the bottom of Figure~11 should
be reconsidered in the light of Theorem~\ref{n}:
namely, at least when~$n\le7$, the $s$-minimal sets
should not look as ``wild'' as they were depicted,
at least for~$s\in ((1/2)-\epsilon_n,\,1/2)$
(and, when~$n=2$, for any~$s\in(0,1/2)$, recall Theorem~\ref{2}):
we try to take into account this further regularity
property in Figure~12, by extending the picture
of the ``nice'' surface down to an unknown exponent~$s=?$.
\medskip

\begin{center}
\includegraphics[width=4in]{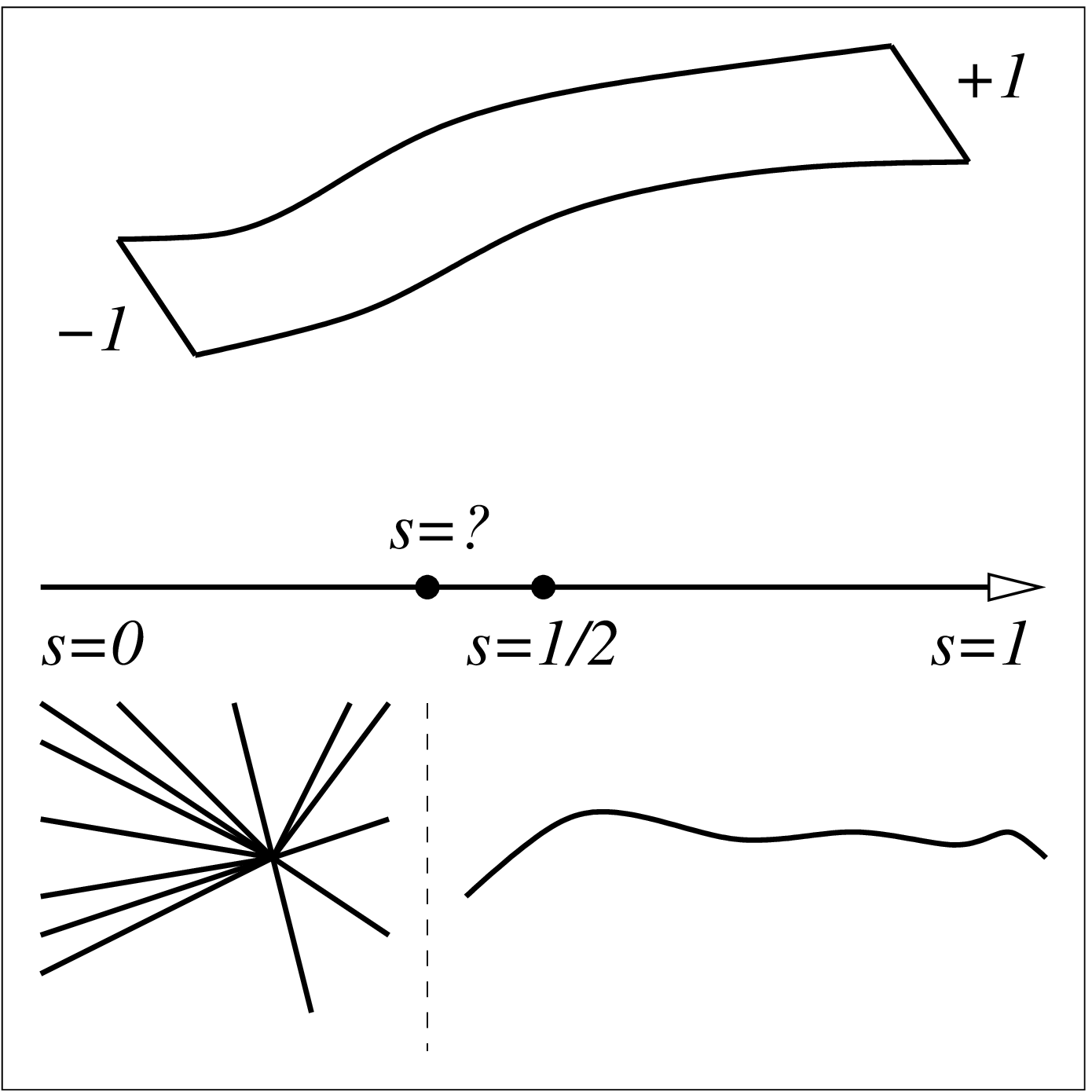}
\end{center}
\nopagebreak\centerline{\footnotesize\em
Figure~12. $\Gamma$-convergence for fractional
phase transitions, taking into account
Theorem~\ref{n}.}
\bigskip

Now we discuss if and how these types of results may
have an influence on the symmetry properties of the solution
of the fractional Allen-Cahn equation.

\subsection{One-dimensional symmetry}

The classical Allen-Cahn equation~\eqref{AC}
is linked to a very famous problem posed by De Giorgi:

\begin{conj}[\cite{DG}]\label{CON}
Let~$u\in C^2(\R^n)\cap L^\infty(\R^n)$ be a solution of
the classical Allen-Cahn equation~\eqref{AC} in the whole of~$\R^n$
and suppose that~$\partial_{x_n} u(x)>0$ for all~$x\in\R^n$.

Then, is it true that $u$ is one-dimensional (i.e., it depends only on
one Euclidean variable up to rotation, and its level sets
are hyperplanes), at least if~$n\le8$?
\end{conj}

Many outstanding mathematicians have given fundamental contributions
to this problem and we cannot do justice here
to all the results obtained and of all the important generalisations
performed (see, e.g.,~\cite{state} for a recent review on the topic):
here we will just mention that Conjecture~\ref{CON} is known to
have an answer in the affirmative when~$n\le3$
and in the negative when~$n>8$.
Also, it is conceivable that Conjecture~\ref{CON}
was inspired by the relation between the phase transitions
and the minimal surfaces (recall Theorem~\ref{sd9fffffff})
and by the rigidity and regularity features of the minimal surfaces.

Of course, a natural question is whether or not
results inspired by Conjecture~\ref{CON} hold true when
the classical Allen-Cahn equation~\eqref{AC}
is replaced by the fractional Allen-Cahn equation~\eqref{AC2}.
This question was addressed in~\cite{sola}
when~$n=2$ and~$s=1/2$, in~\cite{SireV, SireCab} when~$n=2$
and~$s\in(0,1)$, in~\cite{cinti} when~$n=3$ and~$s=1/2$
and in~\cite{cinti2} when~$n=3$ and~$s\in[1/2,1)$. We summarise
these results in the following statement:

\begin{theorem}[\cite{sola, SireV, SireCab, cinti, cinti2}]\label{dg s}
Let either
\begin{equation}\label{p1}
{\mbox{$n=2$ and $s\in(0,1)$}}
\end{equation}
or
\begin{equation}\label{p2}
{\mbox{$n=3$ and $s\in[1/2,1)$.}}
\end{equation}
Let~$u\in C^2(\R^n)\cap L^\infty(\R^n)$ be a solution of
the fractional Allen-Cahn equation~\eqref{AC2} in the whole of~$\R^n$
and suppose that~$\partial_{x_n} u(x)>0$ for all~$x\in\R^n$.
Then $u$ is one-dimensional.
\end{theorem}

Clearly, Theorem~\ref{dg s} leaves many questions open.
For instance, unlike
the classical case, no counterexample
is known in a higher dimension (it is conceivable, but not trivial to 
prove, that the counterexample in~\cite{delpino} works
for~$n>8$ and~$s\in[1/2,1)$; the other ranges of~$n$ and~$s$
seem to be completely unknown).
Furthermore, the ranges in~\eqref{p1} and~\eqref{p2}
are still rather mysterious. One may think that these ranges
are somehow reminiscent of the limit behaviour
of the interfaces, according to Theorem~\ref{d88fdghjww}.
In this spirit, one may suspect that the symmetry result
in Theorem~\ref{dg s} under condition~\eqref{p1}
is a byproduct of the complete regularity theory
of the minimisers of the related perimeter functionals
(i.e., of the classical minimal surfaces when~$s\in[1/2,1)$,
and of the $s$-minimal sets when~$s\in(0,1/2)$, by Theorem~\ref{2}).
Similarly, one may suspect that the threshold~$s=1/2$
of condition~\eqref{p2} is an offspring of the same threshold
that appears in Theorem~\ref{d88fdghjww}, i.e. that the symmetry
properties in Theorem~\ref{dg s} may break down in general 
when~$s\in(0,1/2)$ due to the ``wilderness'' of the limit
$s$-minimal sets. These arguments lead to the feeling that
the threshold $s=1/2$ in~\eqref{p2} is optimal.

On the other hand, some other observations
may lead to an opposite conclusion, that is the feeling that
the threshold $s=1/2$ in~\eqref{p2} may be lowered a bit
(maybe in dimensions~$n=3,\dots,8$). Indeed, the proofs of
Theorem~\ref{dg s} do not explicitly use
any regularity
properties of the (possibly fractional) minimal surfaces
and the threshold $s=1/2$ in~\eqref{p2} does not come from
geometric considerations but from analytical energy estimates.
Moreover, if any relation between symmetry results for phase
transitions and regularity results for (possibly fractional) minimal 
surfaces really holds true, then these
regularity results for $s$-minimisers
hold true up to a threshold~$s=(1/2)-\epsilon_n$,
so it seems conceivable that the
symmetry properties may hold even slightly below~$s=1/2$ at least
when~$n\le8$ (recall Theorem~\ref{n} and Figure~12).
So we believe that any further investigation into
the possible regularity of the $s$-minimal sets and the
possible symmetries of the monotone solutions of the fractional
Allen-Cahn equation would be a 
pleasant challenge for the researchers involved and a
welcome progress for mathematics.

\end{document}